\documentclass[12pt,draft]{amsart}

\newtheorem{theorem}{Theorem}

\newtheorem{thm}[theorem]{Theorem}
 
\newtheorem{proposition}[theorem]{Proposition} 
\newtheorem{cor}[theorem]{Corollary} 
\newtheorem{corollary}[theorem]{Corollary} 
\newtheorem{properties}[theorem]{Properties}

\theoremstyle{definition}

\theoremstyle{remark}

\DeclareMathOperator{\cf}{cf}
\DeclareMathOperator{\pcf}{pcf}

\newcommand{\aaa}{{\mathfrak a}} 
\newcommand{\bbb}{{\mathfrak b}} 
\newcommand{\pcfa}{\pcf \aaa}

\title[decomposable ultrafilters possible cofinalities]
{A connection between decomposable ultrafilters and possible cofinalities. II}

\author[]{Paolo Lipparini} 
\address{Dipartimento di Matematica\\
II Universit\`a di Roma (Sor Vergata)\\
Viale della Ricerca Scientifica\\
I-00133 ROME 
ITALY
}
\email{lipparin@axp.mat.uniroma2.it}
\urladdr{http://www.mat.uniroma2.it/\textasciitilde lipparin}
\thanks{The author has received support from MPI and GNSAGA.
We wish to express our  gratitude to X. Caicedo for stimulating
discussions}



\begin{document} 

\begin{abstract} 
We use 
Shelah's 
 theory of possible cofinalities in order to solve a problem
about ultrafilters.

\begin{thm} \label{abst}
Suppose that $ \lambda $ is a singular cardinal,
$ \lambda ' < \lambda $, and the ultrafilter $D$ is
$ \kappa $-decomposable for all regular cardinals $ \kappa $
with $ \lambda ' < \kappa < \lambda $.  Then $D$ is either
$ \lambda $-decomposable, or $ \lambda ^+$-decomposable.
\end{thm} 

We give applications to topological spaces and to abstract logics
(Corollaries \ref{top}, \ref{log} and Theorem \ref{logics}).
\end{abstract}

\maketitle 

If $F$ is a family of subsets of some set $I$,
 and $ \lambda $ is an infinite cardinal, 
a \emph{$ \lambda $-decomposition} for $F$ 
is a 
function $f: I \to\lambda $
such that whenever $ X \subseteq \lambda $ and 
$ | X|<\lambda $ 
then $ \{ i \in I | f(i) \in X\} \not\in F$. 
The family $F$  is \emph{$\lambda$-decomposable} 
if and only if there is a 
$ \lambda $-decomposition for $F$. 
If $D$ is an ultrafilter (that is, a maximal proper filter) 
let us define the {\it decomposability spectrum} 
$K_D$  of $D$ by $K_D=\{\lambda \geq \omega| D \text{ is }
\lambda\text{-decomposable}\}$.

The question of 
the possible
 values  the spectrum $K_D$ may
take is particularly intriguing.
Even the old problem  from \cite{Si} of characterizing those
$ \mu $ for which there is an ultrafilter $D$ such that 
$K_D=\{\omega ,\mu\} $
is not yet completely solved \cite[p. 1007]{Shr}.

The case when $K_D$ is infinite is even more involved.
\cite{P1} studied the situation in which $ \lambda $ is limit
and $K_D \cap \lambda $ is unbounded in $ \lambda $; he
found some assumptions which imply that $ \lambda \in K_D$.
This is not always the case; if $ \mu $ is
strongly compact and $\cf \lambda < \mu < \lambda $
then there is an ultrafilter $D$ such that 
$K_D \cap \lambda $ is unbounded in $ \lambda $,
and $D$ is not $ \lambda $-decomposable.
If we are in the above situation, we have that necessarily 
$D$ is $ \lambda^+$-decomposable (by \cite[Lemma 3]{So} and the proof of \cite[Proposition 2]{P1}).

The above examples suggest the problem whether 
$K_D \cap \lambda $  unbounded in $ \lambda $
implies that either $ \lambda \in K_D$
or $ \lambda^+ \in K_D$. In general, the problem is still open;
 here we solve it
affirmatively in the particular case when 
there is $ \lambda' < \lambda $ such that $K_D$
contains all regular cardinals in the interval $[\lambda',\lambda )$;
moreover, when $\cf \lambda >\omega $ it is sufficient to assume that 
 $ \{\kappa <\lambda | \kappa^+\in K_D \cap \lambda \}  $ is stationary in $\lambda$.

We briefly review some known results on $K_D$.
If $ \kappa $ is regular and $ \kappa^+ \in K_D$ 
then $ \kappa \in K_D$; and if $ \kappa \in K_D$ is singular, then 
$\cf \kappa \in K_D$. 
Results from \cite{D} imply that
if there is no inner model with a measurable cardinal then
 $K_D$  is always an interval with minimum $ \omega  $.
On the other hand, it is trivial that 
$K_D=\{\mu\} $
if and only if $\mu $ is either $\omega $ or a measurable cardinal.
Further comments and constraints on $K_D$ are given in \cite{L2, L1}. 
Apparently, the problem of determining which sets of cardinals
can be represented as
$ K_F=\{\lambda \geq \omega| F \text{ is }
\lambda\text{-decomposable}\}$ for a filter $F$ has not been studied.  

If $ (\lambda _j) _{j\in J}$ are regular cardinals, the 
\emph{cofinality} $\cf \prod_{j \in J} \lambda_j$
of the product $\prod_{j \in J} \lambda_j$ is the smallest
cardinality of a set
$ G \subseteq  \prod_{j \in J} \lambda_j$
having the property  that for every $f \in \prod_{j \in J} \lambda_j$
there is $g \in G$ such that $f(j) \leq g(j)$ for all $j \in J$.

We shall state our results in a quite general form, involving arbitrary 
filters,
rather than ultrafilters. In what follows, the reader interested in ultrafilters only
can always assume that $F$ is an ultrafilter.


\begin{proposition} \label{prod}    
If $ (\lambda _j) _{j\in J}$ are infinite regular cardinals,
$ \mu= \cf \prod_{j \in J} \lambda_j$
and the filter $F$ is $ \lambda _j$-decomposable
for all $j \in J$, then $F$ is 
$ \mu'$-\hspace{0 pt}decomposable for some $ \mu'$ with
$\sup_{j \in J} \lambda _j \leq \mu' \leq \mu$.
\end{proposition}   

\begin{proof}
Let $F$ be over $I$, and
let $(g_ \alpha) _{ \alpha \in \mu }  $ witness $ \mu= \cf \prod_{j \in J} \lambda_j$. 
For every $j \in J$ let $f(j, -): I \to \lambda _j$ be a $ \lambda _j$
decomposition for $F$. For any fixed $i \in I$,
$f(-, i) \in \prod_{j \in J} \lambda_j$, thus there is
$ \alpha (i) \in \mu$ such that $f(j,i) \leq g _{ \alpha (i)} (j) $
for all $j \in J$.

Let 
$X$ be a subset of $\mu$ with minimal cardinality with
respect to the property that $ Y=\{i \in I | \alpha (i) \in X\} \in F$. Let $\mu'=|X|$.
Thus, whenever $X' \subseteq \mu$ and $|X'|< \mu'$, we have
$ Y'=\{i \in I | \alpha (i) \in X'\} \not\in F$.
Define $h(i)= \alpha (i)$ for $i \in Y$, and 
$h(i)=0$ for $i \not\in Y$. Thus, $h: I \to X \cup \{0\} $.

If $|X'| < \mu'$ then 
$\{i \in I | h(i) \in X'\} \subseteq Y' \cup (I \setminus Y)\not \in F$
(otherwise, since $F$ is a filter, 
$ Y' \supseteq Y \cap Y' =Y \cap (Y' \cup (I \setminus Y)) \in F$, contradiction).
This shows that,
modulo a bijection from $X\cup \{ 0 \}$ onto $\mu'$,
$h$ is a $\mu'$-decomposition for $F$. Trivially, $\mu' \leq \mu$.

Hence, it remains to show that $\sup_{j \in J} \lambda _j \leq \mu'$.
Suppose to the contrary that $\mu' < \lambda _{\bar{j}}$ for some $\bar{j} \in J$.  
Then 
$ |\{ g_{\alpha(i)} (\bar{j}) | i \in Y \} |
\leq |\{ \alpha(i)  | \alpha (i) \in X \} |
 \leq |X| = \mu' < \lambda _{\bar {j}} $. Since $\lambda _{\bar {j}} $ is regular,
we have that
$\beta=\sup _{i \in Y} g _{\alpha(i)} (\bar{j})< \lambda _{\bar{j}}$.
Hence, if $i \in Y$, then $f(\bar{j},i) \leq g _{ \alpha (i)} (\bar{j}) \leq \beta < \lambda _{\bar{j}}$. 
Thus, $|[0, \beta ]|< \lambda _{\bar{j}}$, but 
$ \{ i\in I | f(\bar{j},i) \in [0,\beta ]\} \supseteq Y \in F$, and  
this contradicts the assumption that 
$f(\bar{j}, -)$ is a $ \lambda _{\bar{j}} $
decomposition for $F$.
\end{proof}

Proposition \ref{prod} has not the most general form:
we have results dealing with  the cofinality $\mu$ 
of reduced products $\cf \prod_E \lambda _j$, where $E$
a filter on $J$. We shall not need this more general version here. 

Recall that an ultrafilter $D$ is  $( \mu ,\lambda)$-regular if and only if 
there is a family of $\lambda $ members of $D$
such that the intersection of any $\mu $ members of the family is
empty. We list below the properties of decomposability and regularity
we shall need. Much more is known: see \cite{Lppams, L1} and references there.

\begin{properties} \label{fact}   
(a) Every
$ \lambda $-decomposable ultrafilter is
$\cf \lambda $-\hspace{0 pt}decomposable.
 
(b) Every $\cf \lambda $-decomposable ultrafilter is
$( \lambda , \lambda)$-regular.

(c) If $\mu' \geq \mu$ and $ \lambda' \leq \lambda $ then
every $(\mu,\lambda )$-regular ultrafilter 
is $(\mu',\lambda' )$-regular.

(d) \cite[Theorem 1]{CC} \cite[Theorem 2.1]{KP}
If  $ \lambda $ is singular, $D$ is a
$ \lambda ^+$-decomposable ultrafilter, 
and $D$ is not $\cf \lambda $-decomposable
then $D$ 
is $( \lambda' , \lambda^+ )$-regular
for some $ \lambda ' < \lambda $.

(e) \cite[Corollary 2.4]{Ka}  
If  $ \lambda $ is singular then every
$ \lambda ^+$-decomposable ultrafilter 
is $( \lambda, \lambda^+ )$-regular.

(f) \cite[Corollary 1.4]{Lparch} If  $ \lambda $ is singular then  every
$( \lambda, \lambda)$-regular ultrafilter 
is either $\cf \lambda $-decomposable
or $( \lambda' , \lambda )$-regular
for some $ \lambda ' < \lambda $.

(g) If $ \lambda $ is regular then an ultrafilter is
$ \lambda $-decomposable if and only if it is 
$(\lambda,\lambda )$-regular.
\end{properties}

\begin{theorem}\label{filt}
Suppose that $ \lambda $ is a singular cardinal, 
$F$ is a filter, and either

(a) there is $ \lambda ' < \lambda $ such that  $F$ is
$ \kappa $-decomposable for all regular cardinals $ \kappa $
with $ \lambda ' < \kappa < \lambda $, or  

(b)  $\cf \lambda> \omega $ and 
$ S= \{ \kappa < \lambda | 
 F \text{ is $ \kappa^+ $-decomposable}\} $
is stationary in $ \lambda $.

Then 
$F$ is either
$ \lambda $-decomposable, or $ \lambda ^+$-decomposable.

If $F=D$ is an ultrafilter, then $D$ is 
$( \lambda, \lambda )$-regular.
Moreover, $D$ is
either 
(i) $\lambda $-decomposable,
or
(ii) $( \lambda' , \lambda^+ )$-regular
for some $ \lambda ' < \lambda $,
or
(iii) $\cf \lambda $-decomposable and $( \lambda , \lambda^+)$-regular.
\end{theorem}

\begin{proof}
Recall from \cite{She} that if $\aaa$ is a set of regular cardinals,
then $\pcfa$ is the set of regular cardinals which can be obtained
as $\cf\prod_E \aaa$, for some ultrafilter $E$ on $\aaa$.
If $ \cf \lambda = \nu >\omega $ then by \cite[II, Claim 2.1]{She}
there is a sequence $ (\lambda_ \alpha)_{ \alpha \in \nu }$ 
closed and unbounded in $ \lambda $ and such that, letting 
$\aaa=\{ \lambda^+ _ \alpha | \alpha  \in \nu \}$, we have 
$ \lambda ^+= \max \pcf \aaa $. 
If $ \cf \lambda = \omega $ then 
we have 
$ \lambda ^+= \max \pcf \aaa $
for some countable $ \aaa $ unbounded in $\lambda $
as a consequence of \cite[II, Theorem 1.5]{She}
(since $\aaa$ is countable, any ultrafilter over $\aaa$
is either principal, or extends the dual of the ideal of bounded
subsets of $\aaa$).

Letting $\bbb= \aaa \cap [\lambda ', \lambda )$ in case (a),
and $\bbb= \aaa \cap \{ \kappa^+ | \kappa \in S\} $ in case (b),
we still have 
$\max \pcf \bbb= \lambda ^+$,
because $\bbb$ is unbounded in $ \lambda $, hence
$\max \pcf \bbb \geq \lambda ^+$,
and because $\max \pcf \bbb \leq \max \pcf \aaa = \lambda ^+$,
since $ \bbb \subseteq \aaa$.

Assume, without loss of generality, that $ \lambda' >(\cf \lambda)^+  $
in (a), and that $\inf S > (\cf \lambda)^+ $  in (b).
Since $|\bbb| \leq |\aaa|= \cf \lambda $, then 
$|\bbb|^+< \min \bbb$, hence,
by \cite[II, Lemma 3.1]{She}, 
$ \lambda ^+=\max \pcf \bbb = \cf \prod_{\kappa \in \bbb} \kappa $.
Then Proposition \ref{prod} 
implies that $F$ is either
$ \lambda $-decomposable, or $ \lambda ^+$-decomposable.

The last statements follow from Properties \ref{fact}(a)-(e). 
\end{proof}

\begin{corollary} \label{solovaygenppy}
If $ \lambda  $ is a singular cardinal
and the ultrafilter $D$ is not $ \cf  \lambda $-decomposable,
then the following conditions are equivalent:

(a) There is $ \lambda ' < \lambda $ such that  $D$ is
$ \kappa $-decomposable for all regular cardinals $ \kappa $
with $ \lambda ' < \kappa < \lambda $.

(a$'$) (Only in case $\cf \lambda > \omega $) $\{ \kappa < \lambda | 
 F^+ \text{ is $ \kappa^+ $-decomposable}\} $
is stationary in $ \lambda $.

(b) $D$ is $\lambda^+$-decomposable.

(c) There is $ \lambda' < \lambda $ such that 
$D$ is $(\lambda',\lambda^+ )$-regular.

(d) $D$ is $(\lambda ,\lambda )$-regular.

(e) There is $ \lambda' < \lambda $ such that 
$D$ is $(\lambda',\lambda)$-regular.

(f) There is $ \lambda' < \lambda $ such that
$D$ is $(\lambda'',\lambda'')$-regular
for every $ \lambda ''$ with $ \lambda' < \lambda ''< \lambda $. 
\end{corollary}

\begin{proof} 
 (a) $\Rightarrow$ (b) and (a$'$) $\Rightarrow$ (b) are immediate from Theorem \ref{filt} and Property
\ref{fact}(a).
 In case $\cf \lambda > \omega $, (a) $\Rightarrow$ (a$'$) is trivial.

(b) $\Rightarrow$ (c) $ \Rightarrow $ (d) $\Rightarrow $ (e) $\Rightarrow $ (f) $ \Rightarrow $ (a) are given, respectively, by Properties \ref{fact}(d)(c)(f)(c)(g). 
\end{proof}

\begin{corollary} \label{limit} 
If $ \lambda $ is a singular cardinal, then 
an ultrafilter is $(\lambda ,\lambda )$-regular
if and only if it is either $ \cf \lambda $-decomposable
or $\lambda^+ $-decomposable. 
\end{corollary}
  
\begin{proof}
Immediate from 
 Corollary \ref{solovaygenppy}(d)$ \Rightarrow $(b)
and
Properties \ref{fact}(b)-(d). \end{proof}

A topological space is $[\mu,\lambda   ]$-{\em compact} if and only if every open cover
by $ \lambda $ many sets has a subcover by $<\mu $ many sets. A family $\mathcal F $
of topological spaces is {\em productively} $[\mu,\lambda ]$-compact if and only if every
(Tychonoff) product of members of $\mathcal F$ is
$[\mu,\lambda  ]$-compact.

\begin{cor}
\label{top}
If $ \lambda $ is a singular cardinal, then 
a family of
topological spaces is productively $[\lambda  ,\lambda   ]$-compact
if and only if it is either productively $[ \cf \lambda,  \cf \lambda ]$-compact
or productively $[  \lambda ^+,\lambda^+ ]$-compact. 
\end{cor}

\begin{proof}  
Immediate from Corollary \ref{limit}, Property \ref{fact}(g) and \cite[Theorem 3]{Lptop} (see also \cite{Ca}).
  \end{proof}

Henceforth, by a \emph{logic}, we mean a \emph{regular  logic}  in the sense
of
\cite{Eb}. Typical examples of regular logics are infinitary logics,
 or extensions
of
first-order logic obtained by adding new quantifiers;
e. g., cardinality quantifiers asserting
``there are at least $\omega_\alpha $ $x$'s such that \dots".

A logic $L$ is $[\lambda,\mu ]$-{\em compact}  if and only if for every pair of
sets $\Gamma$ and $\Sigma$ of sentences of $L$, if $|\Sigma|\leq\lambda $
and if $\Gamma \cup \Sigma'$ has a model for every $\Sigma'
\subseteq \Sigma$ with $|\Sigma| < \mu $, then $\Gamma \cup \Sigma$ has a
model
(see \cite{Ma} for some history and further comments).

\begin{cor}
\label{log}
If $ \lambda $ is a singular cardinal, then 
a logic is $[\lambda  ,\lambda   ]$-compact
if and only if it is either $[ \cf \lambda,  \cf \lambda ]$-compact
or  $[  \lambda ^+,\lambda^+ ]$-compact. 
\end{cor}

\begin{proof}  
Immediate from Corollary \ref{limit}, Property \ref{fact}(g) and \cite[Theorem 1.4.4]{Ma} (notice that in \cite{Ma} in the definition of $( \lambda , \mu)$-regularity 
for an ultrafilter 
the order of $\mu$ and $\lambda $ is reversed).
  \end{proof}

\begin{theorem} \label{logics}
Suppose that
$(\lambda_i)_{i \in I}$ and
$(\mu_j)_{j \in J}$
are sets of infinite cardinals. Then the following are equivalent:

(i) For every $i\in I$ there is a $(\lambda_i,\lambda_i)$-regular ultrafilter
which for no $j \in J$ is $(\mu_j,\mu_j)$-regular.

(ii) There is a logic which is $[\lambda_i,\lambda_i]$-compact for every $i\in I$, and
which for no $j \in J$ is $[\mu_j,\mu_j]$-compact.

(iii) For every $i\in I$ there is a $[\lambda_i,\lambda_i]$-compact logic
which for no $j \in J$ is $[\mu_j,\mu_j]$-compact.

The logics in (ii) and (iii) can be chosen to be generated by at most
$2 \cdot |J|$ cardinality quantifiers.
\end{theorem}

\begin{proof} 
In the case when all the $\mu_j$'s are regular, the Theorem is proved
in \cite[Theorem 4.1]{Lparch}. The general case follows from
the above particular case, by applying
Corollaries \ref{limit}  and \ref{log}.  
\end{proof}

\end{document}